\documentclass{amsart}

\newcommand{\ws}{\hspace{4pt}}
\newtheorem{theorem}{Theorem}

\newtheorem{proposition}{Proposition}
\newtheorem{cor}{Corollary}
\newtheorem{lemma}{Lemma}
\newtheorem{obs}{O}

\usepackage[]{amssymb,amsopn}
\usepackage[]{amsmath}
\usepackage[]{eucal}
\usepackage[]{latexsym}

\begin{document}

\title[]{Note on the Equilibrium Measures of Julia sets of Exceptional Jacobi Polynomials}
\author{\'A. P. Horv\'ath }

\thanks{Supported by the NKFIH-OTKA Grant K128922.}

\begin{abstract} We prove that similarly to the standard case, the equilibrium measure of Julia sets of exceptional Jacobi polynomials tends to the equilibrium measure of the interval of orthogonality in weak-star sense.

\end{abstract}
\maketitle

\section{Introduction, Result}

There are several different approaches to define a sequence of measures which tends (at least in weak-star sense) to the equilibrium measure of a compact set $K$ of positive capacity on the complex plane (or on the real line). By orthogonal polynomials with respect to a measure supported on $K$ (with finite moments) it is proved in rather general circumstances, that for instance the normalized counting measure based on the zeros of the orthogonal polynomials in question is an appropriate sequence, see e.g. \cite{ka} and the references therein. When the sequence of densities are the weighted reciprocal of the Christoffel functions, the measures tends to the equilibrium measure of the interval of orthogonality again, see e.g. \cite{mnt}, \cite{ha1}. It is also shown that the normalized counting measure based on the eigenvalues of the truncated multiplication operator has the same limit again, see e.g.  \cite{st}, \cite{s}. Moreover the normalized counting measure based on the zeros of the average characteristic polynomials tend to the equilibrium measure as well, see e.g. \cite{sz}, \cite{ha1}, \cite{ha2}. Dynamical properties of sequences of orthonormal polynomials given by a Borel probability measure supported on a non-polar compact subset of the complex plane were investigated in \cite{chpp} and \cite{pu}.

The results on the different zero counting measures and Christoffel functions are extended to exceptional Jacobi polynomials, see \cite{bo}, \cite{h1}, \cite{h}. The aim of this note is to derive the corresponding result with respect to the equilibrium measure of Julia sets of exceptional Jacobi polynomials. Despite the cited results are rather general, they cannot be applied to exceptional polynomial sequences, since exceptional polynomials are different from the standard or from the classical orthogonal polynomials.  Exceptional orthogonal polynomials are complete systems of polynomials with respect to a positive measure, but the exceptional families have finite codimension in the space of polynomials. Similarly to the classical ones exceptional polynomials are eigenfunctions of Sturm-Liouville-type differential operators but unlike the classical cases, the coefficients of these operators are rational functions. Exceptional orthogonal polynomials also possess a Bochner-type characterization as each family can be derived from one of the classical families applying finitely many Darboux transformations, see \cite{ggm}. Below we deal with exceptional Jacobi polynomials defined by one Darboux transformation and of arbitrary (finite) codimension.\\

Let $w^{(\alpha,\beta)}(x)=(1-x)^\alpha(1+x)^\beta$ and $\left\{\hat{P}_{n}\right\}_{n=0}^{\infty}$ be the system of exceptional Jacobi polynomials orthonormal on $[-1,1]$ with respect to the weight $W=c_0\frac{w^{(\alpha+\varepsilon_1,\beta+\varepsilon_2)}}{\tilde{b}^2}$ ($\varepsilon_i=\pm 1$, $i=1,2$) and generated by one Darboux transformation from the original orthonormal Jacobi polynomials $\{p_n^{(\alpha,\beta)}\}$, and of arbitrary codimension. (For precise definitions see the next section.) From the different results mentioned above to compare with Theorem \ref{T1} below subsequently we use only the next one. We state here as we need it below.

\medskip

\noindent {\bf Theorem A.} {\it  With the notation above if $\alpha, \beta \ge -\frac{1}{2}$ and if $n$ is large enough, $P_n^{[1]}$ has $m$ exceptional zeros in $\mathbb{C}\setminus [-1,1]$, $n$ simple regular zeros in $(-1,1)$. The exceptional zeros tend to the zeros of $\tilde{b}$.\\
For regular zeros we have
\begin{equation}\label{tmte}\tilde{\mu}_n \to \mu_e, \end{equation}
where
\begin{equation}\label{muhu} \tilde{\mu}_n=\frac{1}{n}\sum_{k=1}^{n}\delta_{x_{k,r}}\end{equation}
and $x_{k,r}$, $k=1, \dots , n$ are the regular zeros of $\hat{P}_{n}$. The convergence is meant in weak-star sense.}

\medskip

The previous theorem in several special cases is proved in \cite{gumm} and in a more general case in \cite{bo}.

From a different perspective, that is studying the dynamics related to exceptional Jacobi system, we give a sequence of measures different from the previous ones with the same weak-star limit i.e. the equilibrium measure of $[-1,1]$. The normalization constant $c_0$ is chosen such that $W(x)dx$ be a probability measure on $[-1,1]$.

The basin of attraction for $\infty$ for $\hat{P}_{n}$ is
\begin{equation}\label{O}\Omega_n:=\{z\in \mathbb{C}: \lim_{k\to\infty}\hat{P}_{n}^k(z)=\infty\},\end{equation}
where $\hat{P}_{n}^k=\hat{P}_{n}\circ\hat{P}_{n}\circ\dots \circ\hat{P}_{n}$, that is composition $k$ times. (Similarly $\hat{P}_{n}^{-k}$ denotes the inverse.)
$$K_n:=\mathbb{C}\setminus \Omega_n, \ws \ws \ws J_n:=\partial \Omega_n=\partial K_n$$
are the filled Julia set and the Julia set for $\hat{P}_{n}$, respectively.\\
The main result of this note is

\begin{theorem}\label{T1} Let $\alpha, \beta \ge -\frac{1}{2}$. For $n\in\mathbb{N}$ let $J_n$ be the Julia set of the exceptional orthonormal Jacobi polynomial $\hat{P}_n$, and let $\mu_n$ and $\mu_e$ be the equilibrium measure of $J_n$ and $[-1,1]$, respectively. Then
\begin{equation}\label{1} \mu_n \to \mu_e\end{equation}
in weak-star sense.
\end{theorem}

We remark here that the equilibrium measure of the Julia set of a polynomial is the unique measure of maximal entropy with respect to the polynomial, cf. \cite[Theorem 17.1]{b} and \cite[Theorem 9]{l}. Subsequently we adapt the chain of ideas of \cite{pu} to exceptional case.

\section{Preliminaries}

\subsection{Polynomial Dynamics}
First we note that for each polynomial $p(z)=\gamma z^n+\dots +c$ of degree $n>1$ there is an $R_p$ such that for all $z$ with $|z|>R_p$ $|p(z)|>2|z|$. Thus $\Omega_p$ (cf. \eqref{O}) can be expressed as
\begin{equation}\label{Op}\Omega_p:=\{z\in \mathbb{C}: \lim_{k\to\infty}p^k(z)=\infty\}=\cup_{k\ge 0}p^{-k}(\mathbb{C}\setminus \overline{D(0,R_p)}).\end{equation}
So the filled Julia set of $\hat{P}_n$ is
\begin{equation}\label{kn}K_n=\cap_{k\ge 0}\hat{P}_n^{-k}(\overline{D(0,R_n)}).\end{equation}
As in general, $J_n$ and $K_n$ are compact, completely invariant sets, i.e. $\hat{P}_n^{-1}(J_n)=J_n=\hat{P}_n(J_n)$.

\subsection{Potential Theory}
$\mu$ is a compactly supported Borel probability measure on the complex plane denoted by $\mu\in \mathcal{M}$. The (logarithmic) potential function of $\mu$ is
$$U^{\mu}(z)=\int_{\mathbb{C}}\log\frac{1}{|w-z|}d\mu(w).$$
The energy of $\mu$ is
$$I(\mu)=\int_{\mathbb{C}}U^{\mu}(z)d\mu(z).$$
Let $K\subset \mathbb{C}$ compact.
$$V(K)=\inf\{I(\mu):\mu\in \mathcal{M}, \mathrm{supp}\mu \subset K\}.$$
The capacity of $K$ is
$$\mathrm{cap}K=e^{-V(K)}.$$

\subsection{Exceptional Jacobi Polynomials}
The notion of exceptional orthogonal polynomials is  motivated by problems in quantum mechanics, see e.g.  \cite{gukm0}. It has a rather extended literature, see eg. \cite{ggm} and the references therein. We use the Bochner-type characterization of exceptional polynomials given in \cite{ggm}.

Classical (orthonormal) Jacobi polynomials $\{p_k^{(\alpha,\beta)}\}_{k=0}^\infty$ are eigenfunctions of the second order linear differential operator {with polynomial coefficients}
\begin{equation}\label{de}T[y]=py''+qy',\end{equation}
with eigenvalues $\lambda_n$. $T$ can be decomposed as
\begin{equation}\label{si}T=BA+\tilde{\lambda}, \ws \mbox{with} \ws A[y]=b(y'-wy), \ws B[y]=\hat{b}(y'-\hat{w}y),\end{equation}
where all the coefficients are rational functions. Then the exceptional polynomials are the eigenfunctions of $\hat{T}$, that is the partner operator of $T$, which is
\begin{equation}\label{ka}\hat{T}[y]=(AB+\tilde{\lambda})[y]=py''+\hat{q}y'+\hat{r}y,\end{equation}
where the coefficients of $\hat{T}$ are rational functions cf. \cite[Propositions 3.5 and 3.6]{ggm}.
\eqref{si} and \eqref{ka} imply that
\begin{equation}\label{sk}\hat{T}Ap_n^{\alpha,\beta}=\lambda_nAp_n^{(\alpha,\beta)},\end{equation}
so exceptional (Jacobi) polynomials can be obtained from the classical ones by application of (finitely many) appropriate first order differential operator(s) to the classical (Jacobi) polynomials.

Subsequently we deal with exceptional Jacobi polynomials obtained by one Darboux transformation:
\begin{equation}\label{A} \hat{P}_n:=\frac{1}{\sigma_n}Ap_n^{(\alpha,\beta)}=\frac{1}{\sigma_n}(b\left(p_n^{(\alpha,\beta)}\right)'-bwp_n^{(\alpha,\beta)}),\end{equation}
where $\sigma_n$ is the normalization constant such that $\{\hat{P}_n\}_{n=0}^\infty$ is an orthonormal system on $I=[-1,1]$ with respect to the probability measure
\begin{equation}\label{ew}W:=\frac{c_0pw_0}{b^2}.\end{equation}
To get a polynomial system, $b$ and $bw$ have to be polynomials. Assume that $\deg b\ge \deg bw+1$. In order to the moments of $W$ be finite, $b\neq 0$ on $(-1,1)$. We assume that $b>0$ on $(-1,1)$ and is monic. The normalization constant $c_0$ ensures that the $W$-measure of $I$ is one. We denote by $\tilde{b}$ that part of $b$ which is positive on $[-1,1]$, typically $\tilde{b}(x)=\frac{b(x)}{1-x}$ or  $\tilde{b}(x)=\frac{b(x)}{1+x}$. $\deg\tilde{b}=:m$ which is just the codimension of the exceptional system. The zeros of $\tilde{b}$ is denoted by $Z_{\tilde{b}}$.

The degree of $\hat{P}_n$ is usually greater than $n$. If $n$ is large enough, than $\hat{P}_n$ has $n$ regular and $m$ exceptional zeros, see e.g. \cite{gumm}, \cite{h0}, \cite{h}. Despite the fact that finitely many ones are missing from the sequence of degrees if the set of gaps is admissible, $\left\{\hat{P}_n\right\}_{n=0}^\infty$ is a complete orthogonal system on $L^2_{W,I}$, see \cite{d}.

To prove Theorem \ref{T1} we will compare exceptional and classical Jacobi polynomials. We need some observations.

\medskip

\begin{obs}\label{L1} Denoting by $\gamma_{n,e}$ the leading coefficient of $\hat{P}_n$, we have
\begin{equation}\label{l1} \lim_{n\to \infty}\gamma_{n,e}^{\frac{1}{n}}=2. \end{equation}
\end{obs}

By \eqref{A} $\sigma_n  \hat{P}_n=b\left(p_n^{(\alpha,\beta)}\right)'-bwp_n^{(\alpha,\beta)}$. Recalling that $b$ is monic and denoting the leading coefficient of $bw$ and $p_n^{(\alpha,\beta)}$ by $B$ and $\gamma_n$, respectively, we have
\begin{equation}\label{eh}\gamma_{n,e}=\frac{\gamma_n(n-\epsilon B)}{\delta_n},\end{equation}
where $\epsilon=0$ or $\epsilon=1$ according to the degree of $b$ and $bw$. By \cite[(4.2)]{h} $\sigma_n=\sqrt{c_0(n(n+\alpha+\beta+1)+\tilde{\lambda})}$ (cf. \eqref{ka}). Thus the result can be derived from the corresponding result with respect to the classical Jacobi polynomials, see e.g. \cite[(4.21.6), (4.3.4)]{sz}.

\medskip

Let
$$\hat{\Pi}_n:=\mathrm{span}\{\hat{P}_0, \dots, \hat{P}_n\}.$$
The next observation is
\begin{obs}There is a fixed $s$ such that
\begin{equation}\label{sp}b^2P\in\hat{\Pi}_{n+s}, \end{equation}
where  $P$ is an arbitrary polynomial of degree $n$, and $n$ is large enough. \end{obs}

Indeed, expressing $P=\sum_{k=0}^na_kp_k^{(\alpha,\beta)}$,
$$\int_Ib^2P\hat{P}_lW=\frac{c_0}{\sigma_l}\int_Ib^2(\sum_{k=0}^na_kp_k^{(\alpha,\beta)})(b\sqrt{l(l+\alpha+\beta+1)} p_{l-1}^{(\alpha+1,\beta+1)}-bwp_l^{(\alpha,\beta)})\frac{pw^{(\alpha,\beta)}}{b^2}$$ $$=k_l\sum_{k=0}^na_k\int_I bp_k^{(\alpha,\beta)}p_{l-1}^{(\alpha+1,\beta+1)}w^{(\alpha+1,\beta+1)}-d_l\sum_{k=0}^na_k\int_Ipbwp_k^{(\alpha,\beta)}p_l^{(\alpha,\beta)}w^{(\alpha,\beta)},$$
where $k_l$, $d_l$ are constants depending on $l$, $alpha$ and $\beta$ and $p(x)=1-x^2$, cf. \eqref{de}.
By orthogonality the first sum is zero if $l>\deg b+n+1$, and the second one is zero if $l>\deg bw+2+n$. The computation is finished by recalling the assumption on degrees, and by denoting $\deg b=s$.

To our purpose the observation above is enough. The exact result is that $\tilde{b}^2P$ can be expressed as a linear combination of exceptional polynomials, see \cite[Lemma 1.1]{d}.

\section{Proof of Theorem \ref{T1}}

\begin{lemma}\label{L2}
\begin{equation}\label{l2} \lim_{n\to \infty}\frac{1}{n}\log |\hat{P}_n(z)|=\log |z+\sqrt{z^2-1}|\end{equation}
locally uniformly on $\mathbb{C}\setminus [-1,1]\setminus Z_{\tilde{b}}$, where that branch of the square root is considered which maps positive numbers to positive numbers.
\end{lemma}

\proof  Recalling that $\deg \tilde{b}=m$ and denoting by $\{x_{i,r}\}_{i=1}^n$, $\{x_{j,e}\}_{j=1}^m$, $\{x_{i}\}_{i=1}^n$ the regular zeros of $\hat{P}_n$, the exceptional zeros of $\hat{P}_n$ and the zeros of $p_n^{(\alpha,\beta)}$, respectively, we have
\begin{equation}\label{77}\frac{1}{n}\log |\hat{P}_n(z)|=\frac{1}{n}\log \gamma_{n,e}+\frac{1}{n}\sum_{i=1}^n\log|z-x_{i,r}|+\frac{1}{n}\sum_{j=1}^m\log|z-x_{j,e}|.\end{equation}
Let $z \in \mathbb{C}\setminus [-1,1]$. Then by \eqref{tmte}
$$ \lim_{n\to \infty}\frac{1}{n}\sum_{i=1}^n\log|z-x_{i,r}|=\int_I\log|z-x|d\mu_e(x)$$ $$= \lim_{n\to \infty}\frac{1}{n}\sum_{i=1}^n\log|z-x_{i}|=\log |z+\sqrt{z^2-1}|-\log 2,$$
where we used the corresponding result for classical Jacobi polynomials, and the last equality fulfils by \cite[Theorem 1]{u}. Note, that on each compact set $K\subset \mathbb{C}\setminus [-1,1]$ the functions $f_n(z)=\frac{1}{n}\sum_{i=1}^n\log|z-x_{i,r}|$ are equicontinuous. Indeed, let $z,w \in K$, $\mathrm{dist}(K,I)=d$ and then
$$|f_n(z)-f_n(w)|\le \frac{1}{n}\sum_{i=1}^n\log\left(1+\frac{|z-w|}{|w-x_{i,r}|}\right)\le \frac{|z-w|}{d}.$$
Since the sequence is pointwise convergent, it is uniformly convergent on $K$ as well.

According to Theorem A again, the third sum of \eqref{77} tends to zero locally uniformly on $\mathbb{C}\setminus [-1,1]\setminus Z_{\tilde{b}}$. Indeed, as $K$ is compact and has a positive distance from $Z_{\tilde{b}}$, if $n$ is large enough, then $|\log|z-x_{j,e}||$ is uniformly bounded on $K$. \\
Comparing the first term of the right-hand side of \eqref{77} to \eqref{l1}, the proof is finished.

\medskip

\begin{cor}\label{C1} There is an $\tilde{R}$, such that $K_n \subset D(0,\tilde{R})$ for all $n\in \mathbb{N}$.
\end{cor}

\proof
As $K_n$ is compact for each $n$, it is enough to show that there is an $\tilde{R}$ such that $|\hat{P}_n(z)|>2|z|$ if $|z|>\tilde{R}$ and $n$ is large enough, cf. \eqref{kn}. The required inequality is ensured by \eqref{l2}.

\medskip

Similarly to \cite[Lemma 2.2]{chpp} one can derive

\begin{cor}\label{P1} There exist $R>1$ and $N \in\mathbb{N}$ such that for all $n>N$
\begin{equation}K_n \subset \hat{P}_n^{-1}(\overline{D(0,R)}) \subset D(0,R). \end{equation}
\end{cor}

\proof
Let $R> \max\{1, \tilde{R}\}$ such that $D(0,R)\supset Z_{\tilde{b}}$. Let $\varepsilon:=\inf_{|z|=R}\log |z+\sqrt{z^2-1}|$, which is positive. By Lemma \ref{L2} $\frac{1}{n}\log |\hat{P}_n(z)|\ge \frac{\varepsilon}{2}$ if $n>N$ and $|z|=R$. We can choose $N$ so large that $\log R < N\frac{\varepsilon}{2}$. Thus the generalized minimum principle implies that
$$ \hat{P}_n(\mathbb{C}\setminus D(0,R))\subset \mathbb{C}\setminus \overline{D(0,R)}, \ws \ws \forall \ws n>N,$$
and so
\begin{equation}\label{inv}\hat{P}_n^{-1}(\overline{D(0,R)})\subset D(0,R).\end{equation}
Comparing \eqref{kn}, Corollary \ref{C1} and \eqref{inv} the proof is complete.

\medskip

\begin{proposition}\label{P2} Let $K\subset \mathbb{C}$ compact, $K\cap [-1,1]=\emptyset$. Then there is an $M\in\mathbb{N}$ (depending on $K$, but independent of $n$) such that  for all $\hat{P}_n$, $\mathrm{deg}\hat{P}_n>0$, and any $w\in K_n$
\begin{equation}\label{p2} \mathrm{card}(\hat{P}_n^{-1}(w)\cap K)<M.\end{equation}
\end{proposition}

\medskip

\proof
In view of Corollary \ref{C1} we can take $|w|<R$ with an arbitrary fixed $R$ rather than $w \in K_n$.
We can also assume, that $M$ is large enough. Let $d=d(K,I)$ be the distance of $K$ and $I$, and define $a:=\frac{1}{\sqrt{1+\frac{d^2}{4}}}$. Let $\|b\|_{\infty, [-1,1]}=:A$, and choose $c:= \frac{1}{\kappa(1+R)A^2}$. Here $\kappa$ is a constant ensured by \eqref{l1} such that $\frac{\gamma_{n+s_0,e}}{\gamma_{n,e}}<\kappa$ with a fixed $s_0$, for all $n$. We assume that $M$ is so large that $a^M<c$.\\
Now suppose indirectly that for all $M$ (large enough) and  $N$ there is an $n>N$ such that $\hat{P}_n(z)=w$ has at least $M$ solutions, $x_1, \dots , x_M$  in $K$.\\
Choosing $y_1, \dots , y_M$ to be the nearest points from $I$ to $x_1, \dots , x_M$, respectively, one can define the rational function
\begin{equation}\label{l3} r(z):=\prod_{j=1}^M\frac{z-y_j}{z-x_j},  \ws \ws \|r\|_{\infty,I}\le a^M < c,\end{equation}
cf. \cite[Lemma I.3.2]{stt} and \cite[Lemma 3.3]{pu}.
By this we define the monic polynomial of degree $\deg \hat{P}_n +2s$ as
$$q(z)=\frac{1}{\gamma_{n,e}}rb^2(\hat{P}_n(z)-w).$$
In wiev of \eqref{l3}
$$\|q\|_{2,W}\le \frac{1}{\gamma_{n,e}}\|r\|_{\infty,[-1,1]}\|b^2\|_{\infty, [-1,1]}\|\hat{P}_n(z)-w\|_{2,W}$$ $$\le \frac{1}{\gamma_{n,e}\kappa (1+R)}\|\hat{P}_n(z)-w\|_{2,W}.$$
In that case, when $\deg\hat{P}_0>0$, let $1=\sum_{k=0}^\infty e_k \hat{P}_k$, in norm. By orthonormality
$$\|\hat{P}_n(z)-w\|_{2,W}=\sqrt{1+|w|^2-2\Re w e_n}\le \sqrt{1+R^2+2Re_n}.$$
Note, that if $\deg\hat{P}_0=0$, the right-hand side above is $1+R^2$.
As $1\in L^2_W$, $e_n$ tends to zero thus if $n$ is large enough $\frac{\|\hat{P}_n(z)-w\|_{2,W}}{1+R}<1$.

According to \eqref{sp} There is an $s_0\le s$ such that
\begin{equation}\label{qq} q=\frac{1}{\gamma_{n+s_0,e}}\hat{P}_{n+s_0}+\sum_{k=0}^{n+s_0-1}a_k\hat{P}_{k}.\end{equation}
and choosing $\kappa$ as above we have
$$\|q\|_{2,W}< \frac{\gamma_{n+s_0,e}}{\kappa\gamma_{n,e}}\frac{1}{\gamma_{n+s_0,e}}1< \frac{1}{\gamma_{n+s_0,e}}\|\hat{P}_{n+s_0}\|_{2,W}.$$
Comparing to \eqref{qq} it is impossible.

\medskip

\medskip

\proof (of Theorem \ref{T1}) In view of Corollary \ref{C1} the supports of equilibrium measures of the Julia sets, $\mu_n$, are uniformly bounded. Furthermore applying the M\"obius transform $w_1=\frac{w}{\gamma_{n,e}^{\frac{1}{n-1}}}$, $z_1=\frac{z}{\gamma_{n,e}^{\frac{1}{n-1}}}$ to $w=\hat{P}_n(z)=\gamma_{n,e}z^n+\dots a_0$, and denoting by $F(w):=\frac{1}{n}\sum_{z, \ws \hat{P}_n(z)=w}f(z)$ for an $f$ measurable, and observing that $\frac{1}{n^k}\sum_{z, \ws z= \hat{P}_n^{-k}(0)} f(z)=\frac{1}{n^{k-1}}\sum_{w, \ws w= \hat{P}_n^{-(k-1)}(0)} F(w)$, \cite[Theorem 16.1]{b} reads as
\begin{equation}\label{int}\int_{\mathbb{C}}f(z)d\mu_n(z)=\int_{\mathbb{C}}\frac{1}{n}\sum_{z, \ws \hat{P}_n(z)=w}f(z)d\mu_n(w).\end{equation}
If $K\subset \mathbb{C}$ such that $K$ is compact and disjoint from $I$, then applying \eqref{int} to the indicator function of $K$ and considering Proposition \ref{P2} one can derive that $\mu_n(K)$ tends to zero. Thus, if a subsequence $\mu_{n_k}$ has a weak-star limit, say $\nu$, its support is contained by $I$.

According to \cite[Lemma 15.1]{b} and \eqref{l1}, lower semicontinuity of the energy implies that $I(\nu)\le \liminf I(\mu_{n_k})=I(\mu_e)$. That is $\nu$ must be the unique equilibrium measure for any subsequence and so the proof is finished.

\subsection{Final Remark}
On one side the main ingredients of the proof of Theorem \ref{T1} are the two facts which ensures the result in the standard case as well: that the orthonormal polynomials in question has regular zero-distribution, that is the normalized counting measures based on the zeros tend to the equilibrium measure of the set of orthogonality, $S$, and $\limsup \gamma_n^{\frac{1}{n}}=\frac{1}{\mathrm{cap} S}$. In many cases these facts are connected, see e.g. \cite[Theorem 2.2.1]{stt}, and can be proved by potential theory. Here we derived the the properties of the exceptional system from the well-known properties of classical Jacobi polynomials, but of course, on the right-hand side of \eqref{l2} stands the Green function of the corresponding domain, and of \eqref{l1} the reciprocal of the capacity of $I$.\\
On the other side the finite codimension implies the necessity of a formula like \eqref{sp}. In addition it is enough to require that the constat 1 can be expand to a convergent series.\\
So the next version of Theorem \ref{T1} could be formulated (with the example of exceptional Jacobi systems)

{\bf Theorem 1V.}
{\it If a sequence of polynomials $\{Q_n\}$ which are orthonormal with respect to a compactly supported probability measure $\mu$ ($\mathrm{cap}S>0$) on the complex plane fulfils that the zero-distribution of the polynomials is regular, the leading coefficients, $\gamma_n$ fulfils that $\limsup_{n\to\infty} \gamma_n^{\frac{1}{n}}=\frac{1}{\mathrm{cap}S}$, 1 can be expand to a convergent series, and there is a fixed $s$ $C$ such that for any polynomial $P$ of degree $n$ there is another polynomial $g_P$  such that $\|g_P\|\le C$ and $g_PP\in \Pi_{n+s}=\mathrm{span}\{Q_0, \dots, Q_{n+s}\}$, then the equilibrium measures of the Julia sets of the polynomials tends to the equilibrium measure of $S$.}

\medskip

\noindent \small{Department of Analysis, \newline
Budapest University of Technology and Economics}\newline

\small{ g.horvath.agota@renyi.hu}
\end{document}